\RequirePackage[l2tabu, orthodox]{nag}
\documentclass[10pt, draft]{article}

\usepackage{amsmath,amssymb,amsbsy,amsfonts,amsthm,latexsym,amsopn,amstext,
                            amsxtra,euscript,amscd}
\usepackage{booktabs}
\usepackage{enumitem}
\usepackage{tikz}
\usepackage[norefs,nocites]{refcheck}

\numberwithin{equation}{section}
\newtheorem{theorem}{Theorem}[section]
\newtheorem{lemma}[theorem]{Lemma}

\newtheorem{lem}[theorem]{Lemma}

\DeclareMathOperator{\lcm}{lcm}

\def\\{\vspace{\baselineskip}}
\def\({\left(}
\def\){\right)}
\def\[{\left[}
\def\]{\right]}
\def\<{\langle}
\def\>{\rangle}
\def\floor#1{\left\lfloor#1\right\rfloor}

\def\notdivides{\mathrel{\kern-3pt\not\!\kern3.5pt\bigm|}}

\begin{document}


\title{\Large \textbf{A summation of the number of distinct prime divisors of the lcm}}
\author{%
\scshape{RANDELL HEYMAN}\\
{School of Mathematics and Statistics, University of New South Wales} \\
{Sydney, Australia}\\
\text{randell@unsw.edu.au}
}

\maketitle


\begin{abstract}
Let $x$ be a positive integer. We give an asymptotic result for $\omega(\lcm(m,n))$ summed over all positive integers $m$ and $n$ with $mn \le x$. This answers an open question posed in a recent paper.

NOTE: This paper has now been incorporated in R. Heyman and \newline L. T\'oth, `On certain sums of arithmetic functions involving the gcd and lcm of two positive integers', 
\end{abstract}

\newline


\section{Introduction}
Let $x,m$ and $n$ be positive integers throughout. As usual we use $(*.*)$  to notate the greatest common divisor and $[*.*]$ to notate the least common multiple. Results for many summations of the form $\sum f((m,n))$ and $\sum f([m,n])$, where $f$ is an arithmetic function, have been given in the literature. A recent brief summary can be found in \cite{Hey}. In that paper the authors posed an open question to seek an asymptotic formula for
$$S_{\omega}(x):=\sum_{mn \le x}\omega([m,n]).$$
This question can be answered in (at least) two ways.
The first way recognises that
 $$S_{\omega}(x)=\sum_{n \le x} \tau(n) \omega(n),$$
where, as usual, the tau function is given by
$$\tau(n)=\sum_{d|n}1.$$
This, using a known asymptotic result, gives the following:

\begin{theorem}
\label{thm:first}
We have
$$S_\omega(x)= 2x \log x \log \log x +Bx \log x +O(x),$$
where
$$B = 2\(\gamma -1+\sum_p\(\log\(1-\frac{1}{p}\)+\frac{1}{p}-\frac{1}{2p^2}\)\).$$
\end{theorem}
In this theorem and throughout $\gamma$ is the Euler-Mascheroni constant and $\sum_p$ is the summation over all primes.
The second way, based on an elementary relationship between gcd and lcm summations involving the omega function, gives an inferior error term as follows:
\begin{theorem}
\label{thm:second}
We have
$$S_\omega(x)= 2x \log x \log \log x +Bx \log x +O(\log \log x).$$
\end{theorem}

\section{Preparatory lemmas}
We require the following lemmas:
\begin{lemma}
\label{lem:omegas}
Let $k$ be a positive integer. Then
\begin{align*}
\sum_{mn = k}\omega([m,n])&=\sum_{mn = k} \omega(m)+ \sum_{mn=k}\omega(n) -\sum_{mn = k} \omega((m,n))\\
&=2\sum_{d|k}\omega(d)-\sum_{mn = k} \omega((m,n)).
\end{align*}
\end{lemma}
\begin{proof}
If $p$ is a prime divisor of $[m,n]$ then $p$ is also of $m,n$ and $[m,n]$. Conversely, if $p$ is not a prime divisor of $[m,n]$ then $p$ is not a prime divisor of $m,n$ nor $[m,n]$. This proves the first equation. The second equation follows using
$$\sum_{mn = x} \omega(m)= \sum_{mn=k}\omega(n)=\sum_{d|k}\omega(d).$$
\end{proof}
The following is given in \cite{Hey}:
\begin{lemma}
\label{lem:Hey}
We have
\begin{align}
\label{eq:heytoth}
\sum_{mn \le x} \omega((m,n))=Cx\log x - Dx +O\(x^{1/2}\),
\end{align}
where
$$C=\sum_{p} \frac1{p^2},\,\,D=(2\gamma-1)\sum_{p} \frac1{p^2}-2\sum_{p}\frac{\log p}{p^2}.$$
\end{lemma}
The following was first proven in \cite{Saf} using Dirichlet's hyperbola method to prove a more general result:
\begin{lem}
\label{lem:Har}
For any positive integer $h$ we have
$$\sum_{n \le x} \omega(n) = x \log \log x +Ax+ \sum_{j \le h}\frac{a_jx}{\log^jx}+O\(\frac{x}{\log^{h+1} x}\),$$
where $$A = \gamma +\sum_p(\log (1 - 1/p) +1/p), \,\, a_j=-\int_1^\infty \frac{\{t\}}{t^2}(\log t)^{j-1}\, dt.$$
\end{lem}
\section{Proof of Theorem \ref{thm:first}}
We have
\begin{align}
\label{eq:tau proof}
\sum_{mn \le x}\omega([m,n])&=\sum_{k \le x} \sum_{mn=k} \omega([m,n])\notag\\
&=\sum_{k \le x} \sum_{m|k}\omega([m, k/m])
\end{align}
Expressing $k$ as the product of distinct primes we have
$$k=\prod_{i \le \omega(k)}p_i^{\alpha_i},\,\, \alpha_i>0.$$
Since $m$ divides $k$ we have
$$k/m=\prod_{i \le \omega(k)}p_i^{\beta_i},\,\, 0 \le \beta_i \le \alpha_i.$$
It follows that
$$[m,k/m]=\prod_{i \le \omega(k)}p_i^{\max\{\beta_i, \alpha_i-\beta_i\}}.$$
Since $ \max\{\beta_i, \alpha_i-\beta_i\}>0$ for all $i$ it follows that
$$\omega([m, k/m])=\omega(k).$$
Substituting this into \eqref{eq:tau proof} we have
$$\sum_{mn \le x}\omega([m,n])=\sum_{k \le x} \sum_{m|k}\omega(k)=\sum_{k \le x} \omega(k)\sum_{m|k}1=\sum_{k \le x}\tau(k)\omega(k).$$
We then use the error term in the asymptotic result in \cite[Theorem 9]{DeK}, which completes the proof of the theorem.

\section{Proof of Theorem \ref{thm:second}}
Using Lemma \ref{lem:omegas} we have
\begin{align}
\label{eq:main thm}
S_\omega(x)&=\sum_{ k \le x}\sum_{mn=k} \omega([m,n])\notag\\
&=2 \sum_{k \le x}\sum_{d|k}\omega(d)-\sum_{mn \le x} \omega((m,n)).
\end{align}
Next,
\begin{align*}
2\sum_{k \le x}\sum_{d|k}\omega(d)&=2\sum_{k \le x}\sum _{d \le x/k} \omega(d)\\
&=2(S_1+S_2),
\end{align*}
where
$$S_1=\sum_{\floor{x/2}<k\le x}\sum_{d\le x/k}\omega(d)$$
and
$$S_2=\sum_{k \le \floor{x/2}} \sum_{d\le x/k}\omega(d).$$
For $S_1$, the summation limit that $k>\floor{x/2}$ implies that $d=1$ in all cases. So every $\omega(d)=0$ and thus $S_1=0$.
So, using Lemma \ref{lem:Har}, we have
\begin{align}
\label{eq:2wd}
2\sum_{k \le x}\sum_{d|k}\omega(d)&=2S_2\notag\\
&=2x(S_{2,1}+S_{2,2}+S_{2,3}+S_{2,4}),
\end{align}
where
$$S_{2,1}=\sum_{k \le \floor{x/2}}\frac{1}{k} \log \log \(\frac{x}{k}\),$$
$$S_{2,2}=\sum_{k \le \floor{x/2}}\frac{A}{k},$$
$$S_{2,3}=\sum_{k \le \floor{x/2}}\sum_{j \le h}\frac{a_j}{k\log^j (x/k)},$$
$$S_{2,4}=O\(\sum_{k \le \floor{x/2}}\frac{1}{k \log^{h+1} (x/k)}\).$$
By comparing $S_{2,1}$ to the appropriate integrals, we have
\begin{align*}
S_{2,1}&\le  \log \log x +\int_1^{x/2}\frac{\log \log (x/k)}{k}\,dk\\
&=\log x \log \log x - \log  x +\log \log x +\log 2 - \log 2 \log \log 2.
\end{align*}
and
\begin{align*}
S_{2,1}&\ge \int_0^{\floor{x/2}-1}\frac{\log \log (x/(k+1))}{k+1}\,dk \ge \int_1^{x/2}\frac{\log \log (x/k)}{k}\,dk\\
&= \log x \log \log x- \log x +\log 2-\log 2 \log \log 2.
\end{align*}
We conclude that
\begin{align}
\label{eq:S21}
S_{2,1} =  \log x \log \log x-\log x +O\( \log \log x\).
\end{align}
Next, using , for example, \cite [Theorem 3.2]{Apo},
\begin{align}
\label{eq:S22}
S_{2,2}&=A\sum_{k \le \floor{x/2}}\frac1{k}=A\log x+O(1).
\end{align}
Turning to $S_{2,3}$. We have, 
by comparison to the appropriate integrals,
\begin{align}
\label{eq:S23}
S_{2,3}&=\sum_{j \le h}a_j\sum_{k \le \floor{x/2}}\frac{1}{k\log^j (x/k)}\notag\\
&\le h \max\{a_1,\cdots,a_h\}\sum_{k \le \floor{x/2}}\frac{1}{\log x}\notag\\
&=O(\log \log x).
\end{align}
Finally,
\begin{align}
\label{eq:S24}
S_{2,4}=O\(\sum_{k \le \floor{x/2}}\frac{1}{k\log^{h+1} (x/k)}\)=O\(\sum_{k \le x/2}\frac{1}{k\log^{h+1} (x/k)}\)=O\(1\).
\end{align}
Substituting \eqref{eq:S24}, \eqref{eq:S23}, \eqref{eq:S22} and \eqref{eq:S21} into \eqref{eq:2wd} we have
\begin{align*}
2\sum_{k \le x} \sum_{d|k} \omega(d)&=2x\log x \log \log x +(2A-2)x \log x+O(x \log \log x).
\end{align*}
Substituting this result and \eqref{eq:heytoth} into \eqref{eq:main thm} we obtain
$$S_\omega(x)= 2x \log x \log \log x +(2A-C-2)x \log x +O(x\log \log x),$$
which concludes the proof.

\section{Acknowledgement}
The author thanks Olivier Bordell\`es for bringing to the author's attention the asymptotic result for the partial sum of the omega function (Lemma \ref{lem:Har}).

\makeatletter
\renewcommand{\@biblabel}[1]{[#1]\hfill}
\makeatother

\end{document}